\documentclass[10pt]{amsart}
\usepackage{amsmath}
\usepackage{amsmath}
\usepackage{amsfonts}
\usepackage{amssymb}

\title{Virtual manifolds and localization}
\author{Bohui Chen}
\address{Department of Mathematics, Sichuan University,
        Chengdu,610064, China}
        \email{xfzhu1@yahoo.com}
\author{Gang Tian}
\address{Department of Mathematics, Princeton University, tian@}
\email{tian@math.princeton.edu}
\usepackage{amsmath,amsthm}
\begin{document}
\maketitle

\newtheorem{theorem}{Theorem}[section]
\newtheorem{assertion}[theorem]{Assertion}
\newtheorem{claim}[theorem]{Claim}
\newtheorem{conjecture}[theorem]{Conjecture}
\newtheorem{corollary}[theorem]{Corollary}
\newtheorem{defn}[theorem]{Definition}
\newtheorem{example}[theorem]{Example}
\newtheorem{figger}[theorem]{Figure}
\newtheorem{lemma}[theorem]{Lemma}
\newtheorem{prop}[theorem]{Proposition}
\newtheorem{remark}[theorem]{Remark}
\newtheorem{assumption}[theorem]{Assumption}

\def \s{stabilization}
\def \mcb{\mathcal{B}}
\def \mcf{\mathcal{F}}
\def \mc{\mathcal}
\def \v{\vskip 0.1in}
\def \n{\noindent}
\def \mcu{\mathcal{U}}
\def \mco{\mathcal{O}}
\def \mfk{\mathfrak}

\def \inv{^{-1}}

\abstract In this paper, we explore the virtual technique that is
very useful in studying moduli problem from differential geometric
point of view. We introduce a class of new objects "virtual
manifolds/orbifolds", on which we develop the integration theory.
In particular, the virtual localization formula is obtained.
\endabstract

\section{Introductions}\label{sect_1}

In this paper, we introduce a class of new objects, which we call
them ``virtual manifolds/orbifolds''.  As the terminology
suggests, it is a generalization of manifold/orbifold. One of the
main themes of this paper is to show that one can do most analysis
on those objects as one does on usual manifolds, particularly, in
we develop a modified integration theory on and show an analogue
of the deRham theory for virtual manifolds/orbifolds. Furthermore,
we study $G$-actions on virtual manifolds. We introduce a notion
of $G$-virtual manifolds/orbifolds and develop a $G$-equivariant
(integration) theory on them. One of the main results in this
paper is the Atiyah-Bott type localization formula on $G$-virtual
manifolds when $G$ is abelian. We call such a formula the "virtual
localization formula" (Theorem 6.8).

Virtual manifolds/orbifolds provide a natural frame to study
certain type of singular spaces that come from moduli problems in
geometry. By "a moduli problem", we mean the construction of
invariants on moduli spaces that are associated to Fredholm
systems (cf. \S\ref{sect_5}). There are many famous moduli
problems of this sort, e.g, the moduli space of anti-self-dual
instantons in defining the Donaldson invariants, the moduli space
in constructing the Seiberg-Witten invariants, the moduli space of
$J$-holomorphic maps (from Riemann surfaces to symplectic
manifolds) constructing the Gromov-Witten invariants. Let us take
the moduli space of $J$-holomorphic maps as an example. The
Gromov-Witten invariants were first constructed for semi-positive
symplectic manifolds (cf. \cite{Ru}. \cite{Ru-Ti}, \cite{Ru-Ti2}).
Since the involved moduli spaces may be {\em singular}, one needs
to introduce the technique of ``virtual cycles'' in order to
construct the Gromov-Witten invariants for general symplectic
manifolds. In around 1996, several groups of people gave different
constructions of virtual cycles. These groups include
Fukaya-Ono(\cite{FO}), Li-Tian(\cite{LT},\cite{LT2}),
Liu-Tian(\cite{Liu}), Ruan{\cite{R}} and etc.. In this paper, we
explain that for any Fredholm system, we are able to construct a
virtual manifold/orbifold associated to its corresponding moduli
problem. The invariants then can be defined via the integration on
this virtual object. Such a general construction can be applied to
the Gromov-Witten theory to
 get  these "virtual cocycles" in the symplectic
case  and therefore get the Gromov-Witten invariants. These are
done in \cite{CL}.  This approach by using integration follows the
one used by Ruan in his construction of the Gromov-Witten
invariants for general symplectic manifolds {\cite{R}). In some
sense, one may also treat the theory of "virtual manifolds" as a
dual to Fukaya-Ono's construction of Kuranishi structures or
Li-Tian's construction of weakly smooth structures. Our
construction can be also carried out for ``weakly'' Fredholm
systems which are more general than Fredholm ones and require less
smoothness. The problem of constructing the Gromov-Witten
invariants is on of such systems.

We then go further to consider Fredholm systems with $G$-actions.
The virtual-manifolds associated to moduli spaces then turn to be
$G$-virtual manifolds. Therefore, we  develop the abelian virtual
localization formula for moduli problems when $G$ is abelian. This
is applied to derive the symplectic virtual localization formula
for Gromov-Witten invariants  in \cite{CL}. We should point out
that such a formula was previously developed in the algebraic
geometry category (\cite{GR}).

\v {\em Acknowledge. }The first author would like to express his
special thanks to A. Li and Y. Ruan for their all time support and
encouragement.

\section{Virtual Orbifolds}\label{sect_2}
In this section, we  introduce a class of new objects "virtual
manifolds(orbifolds)".

\subsection{What is a virtual orbifold supposed to be?}\label{sect_2.1}

An $n$-dimensional manifold/orbifold can be constructed by
patching several pieces of $n$-dimensional manifolds/orbifolds
together-- note that this is not obvious for orbifolds.  From this
point of view, a virtual manifold(orbifold) is obtained by
patching several pieces of {\em possibly different dimensional}
manifolds (orbifolds) together properly.

We take a  simplest example to explain what  we mean by patching.
let $A_1$ and $A_2$ be two manifolds/orbifolds with dimension $n$
and $n+k$ respectively. Let $U_i\subseteq A_i, i=1,2,$ be two open
submanifolds/orbifolds of $A_i$ and suppose that $\pi: U_2\to U_1$
is a rank $k$ (orbifold) vector  bundle. So $U_1$ is identified
with the $0$-section in $U_2$,
 say $U_1'$.
By patching $A_1$ and $A_2$ together, geometrically, we mean the
quotient space
$$
A_1\cup A_2/ (U_1\cong U_1').
$$
Such an object is a virtual manifold/orbifold.

To emphasis the role of the bundle structure $(U_2,U_1,\pi)$, we
introduce a  space $A_1\cup A_2/\sim$, where the equivalent
relation is given by
$$
x\in U_1\sim y\in U_2 \iff \pi(y)=x.
$$
This space is called the virtual space of the virtual
manifold(orbifold) given above.

In the above example, we say that $U_i\subseteq A_i$ are the
overlapping areas in the sense of "patching". To summarize, {\em
two different dimensional manifolds $A_1$ and $A_2$ are patched at
$U_i\subseteq A_i$ which are different up to a vector bundle
structure.}

Technically, the formalism of such objects is not obvious. When
more than two pieces are patched, certain compatibility is needed.
For this purpose, we explain a useful, but rather obvious,
principle of patching in next subsection.

\subsection{A principle of patching}\label{sect_2.2}
Let $N=\{1,\ldots,n\}$ and $\mc N=2^N$ be the set of all subsets
of $N$. Let
$$
\mc X= \{X_I|I\in \mc N\}
$$
be a collection of sets indexed by $\mc N$.
For any $I\subset J$ there exist $ X_{I,J}\subset X_I,
X_{J,I}\subset X_J $ and a surjective map
$$
\phi_{J,I}: X_{J,I}\to X_{I,J}.
$$
Set ${\Phi}=\{\phi_{J,I}|I\subset J\}$. We always assume that $X_\emptyset \not=\emptyset$.
\begin{defn}\label{defn_2.2.1}
A pair $(\mc X,{\Phi})$ is called {\em patchable} if for any
$I,J\in \mc N$ we have
\begin{itemize}
\item[P1.] $X_{I\cup J,I\cap J}= X_{I\cup J,I}\cap X_{I\cup J,J}$;
\item[P2.] $X_{I\cap J,I\cup J}= X_{I\cap J,I}\cap X_{I\cap J,J}$;
\item[P3.] $\phi_{I\cup J,I\cap J}= \phi_{I,I\cap J}\circ
\phi_{I\cup J,I} =\phi_{J,I\cap J}\circ \phi_{I\cup J,J}$;
\item[P4.] $\phi_{I\cup J,I}(X_{I\cup J,I\cap
J})=\phi\inv_{I,I\cap J}(X_{I\cap J,I\cup J})$; \item[P5.]
$\phi_{I\cup J,J}(X_{I\cup J,I\cap J})=\phi\inv_{J,I\cap
J}(X_{I\cap J,I\cup J})$.
\end{itemize}
Set
\begin{eqnarray*}
X_{I,J}&=&\phi_{I\cup J,I}(X_{I\cup J,I\cap J})=\phi\inv_{I,I\cap J}(X_{I\cap J,I\cup J}),\\
X_{J,I}&=&\phi_{I\cup J,J}(X_{I\cup J,I\cap J})=\phi\inv_{J,I\cap
J}(X_{I\cap J,I\cup J}).
\end{eqnarray*}
\end{defn}
In this paper, we always assume that $(\mc X, \Phi)$ is patchable.
We define a relation for points in $\cup X_I$.
\begin{defn}\label{defn_2.2.2}
For $x\in X_I$ and $y\in X_J$ we say that $x\sim y$ if and only if
there exists a $K\subseteq I\cap J$ such that
$$
\phi_{I,K}(x)=\phi_{J,K}(y).
$$
\end{defn}
\n We claim that "$\sim$" is an equivalence relation. This follows
from the next lemma.
\begin{lemma}\label{lemma_2.2.1}
Let $x\in X_I, y\in X_J$ and $z\in X_K$. If $x\sim y$ and $y\sim
z$, then $x\sim z$.
\end{lemma}
{\bf Proof. }By assumptions, we have
\begin{eqnarray*}
&&\phi_{I, K_1}(x)= \phi_{J,K_1}(y);\\
&&\phi_{J,K_2}(y)=\phi_{K,K_2}(z);
\end{eqnarray*}
for some $K_1\subseteq I\cap J$ and $K_2\subseteq J\cap K$. Then
$$
\phi_{I,K_1\cap K_2}(x) =\phi_{J,K_1\cap K_2}(y)=\phi_{K,K_1\cap
K_2}(z).
$$
Clearly, $K_1\cap K_2\subset I\cap K$. Therefore $x\sim z$. q.e.d.

\v We "patch" $X_I$ together and get a set
$$
\mathbf{X}= \bigcup_{I\in \mc N}X_I/\sim.
$$
From a patchable $(\mc X,\Phi)$ to $ {\mathbf X}$ is our so-called
{\em principle of patching} in this paper.

\subsection{Virtual manifolds/orbifolds}\label{sect_2.3}

A virtual manifold is a patchable pair $(\mc X,\Phi)$  with
specified properties.
\begin{defn}\label{defn_2.3.1}
Let $(\mc X,\Phi)$ be a patchable pair. Suppose that
\begin{itemize}
\item $X_I\in \mc X$ are  smooth orbifolds; \item $X_{I,J}$ and
$X_{J,I}$ are  {\em open} suborbifolds in  $X_I$ and $X_J$
respectively; \item $\Phi_{J,I}: X_{J,I}\to X_{I,J}$ is an
orbifold {\em vector bundle}.
\end{itemize}
Then $(\mc X,\Phi)$ is called a {\em virtual orbifold} if for any
$I$ and $J$,
\begin{eqnarray*}
&&\phi_{I,I\cap J}: X_{I,J}\to X_{I\cap J,I\cup J},\\
&&\phi_{J,I\cap J}: X_{J,I}\to X_{I\cap J,I\cup J}
\end{eqnarray*}
are orbifold vector bundles and
\begin{equation}\label{eqn_2.3.1}
X_{I\cup J,I\cap J}= X_{I,J}\times_{X_{I\cap J,I\cup J}} X_{J,I}.
\end{equation}
We call
$$
{\mathbf X}=\bigcup_{I\in \mc N} X_I/\sim
$$
 the {\em virtual space} of $(\mc X,\Phi)$.
We denote the projection map $X_I\to {\mathbf X}$ by $\phi_I$.

Let $d_I$ be the dimension of $X_I$. We call $d_\emptyset $ the
virtual dimension of $(\mc X,\Phi)$.
\end{defn}
For simplicity, from now on, we assume that $(\mc X,\Phi)$ is a
virtual manifold. The discussion is identical for virtual
orbifolds.

A point in ${\mathbf X}$ is an equivalence class, denoted by $[
x]$. Set $[X_I]=\phi_I(X_I)$. Then $\{[X_I]\}$ forms a cover of
$\mathbf X$. For any $[x]\in \mathbf X$, there exist some $X_I$
such that $\phi_I^{-1}([x])$  consists of only
 one single point $x$. Furthermore,
among them, there is a unique  $X_I$ such that $i=|I|$ is
smallest. For $X_I$ and $X_J$ are two such sets, so is $A_{I\cap
J}$. This contradicts to the assumption of smallest.
 We call such an  $X_I$
the {\em support} of $[x]$.

\begin{remark}\label{rmk_2.3.1}
We can define the {\em virtual manifolds with boundary} by a
slight modification of definition \ref{defn_2.3.1}: (1), $X_I$ are
manifolds (possibly with boundaries); (2), if $X_{I,J}$ contains
boundary $\partial X_{I,J}\subseteq \partial X_I$, we require that
$$
\partial X_{J,I}= \Psi_{J,I}^{-1}(\partial X_{I,J}).
$$
Such an object is called a virtual manifold(orbifold) with
boundary. Set
$$
\partial \mc X=
\{\partial X_I|I\in \mc N\}
$$
and
$$
\partial\phi_{I,J}=\phi_{I,J}|_{\partial X_{I,J}}.
$$
Let ${\partial\Phi}=\{\partial \phi_{I,J}\}$. Then $(\partial\mc
X,\partial\Phi)$ also forms a virtual manifold. We call it the
{\em boundary} of $(\mc X,\Phi)$. The induced virtual space
$\partial{\mathbf X}$ is called the {\em  boundary} of $\mathbf
X$.
\end{remark}
We say $[x]$ is  an {\it interior point} if $x$ is an interior
point in its support. Let $\mathbf X^\circ$ denote the set of
interior points of $\mathbf X$. We see that
$$
\mathbf X^\circ = \mathbf X-\partial\mathbf{ X}.
$$
If $\partial\mathbf{ X}$ is empty, we say that $(\mc{X},\Phi)$,
 is {\it boundary free}. If $\mathbf X$ is compact,
We say that $(\mc{X},\Phi)$ is {\it compact}.

We give  examples of virtual manifolds. If $X$ is a manifold.
Itself is clearly a virtual manifold: let $N=\emptyset$,
$A_\emptyset=X$. However, we can construct a nontrivial virtual
manifold out of $X$. This is explained in the following example.
\begin{example}\label{ex_2.3.1}
 Let $X$ be a manifold. Let $\{U_0,U_1, \ldots,U_n\}$ be an open
cover of $X$. Let $U_i^\circ = \frac{3U_i}{4}, i\geq 1$. Here
$\frac{3U_i}{4}$ just means an open subset whose closure is in
$U_i$. We use $\frac{3}{4}$ to make the notations more suggestive.

Let $N=\{1,\ldots, n\}$ and $I,J,K$ be as before. Define
\begin{eqnarray*}
X_\emptyset &=& U_0 - \bigcup_{i=1}^n U_i^\circ \\
X_I &=& \bigcap_{i\in I} U_i - \bigcup_{j\not\in I}U_j^\circ.
\end{eqnarray*}
Let $\mc X= \{X_I|I\in \mc N\}$. Define
$$
X_{I,J}=X_{J,I}= X_I\cap X_J.
$$
All possible $\psi_{J,I}$ are taken to be identities and let
$\mathbf \Phi=\{\phi_{J,I}\}$. Then $(\mc X,\mathbf\Phi)$ is a
virtual manifold (cf. Proposition \ref{prop_2.3.1}). Moreover, the
virtual space $\mathbf X$ is $X$.
\end{example}
\begin{prop}\label{prop_15.2.1}
$(\mc X,\mathbf\Phi)$ given in Example \ref{ex_2.3.1} is a virtual
manifold.
\end{prop}
{\bf Proof: } Let $U_i^c= X-U_i^\circ$. By definition,
$$
X_I= \bigcap_{i\in I}U_i\bigcap_{j\not\in I}U_i^c.
$$
If $I\subseteq J$
\begin{equation}\label{eqn_2.3.2}
X_{I,J}=X_{J,I} =\bigcap_{i\in I}U_i\bigcap_{j\not\in J}U_i^c
\bigcap_{k\in J-I}(U_k-U_k^\circ).
\end{equation}
Now for arbitrary $I,J$, with computations
\begin{eqnarray*}
X_{I\cup J,I}&=& \bigcap_{i\in I}U_i\bigcap_{j\not\in I\cup
J}U^c_j
\bigcap_{k\in J-I} (U_k-U_k^\circ);\\
X_{I\cup J,J}&=& \bigcap_{k\in J}U_k\bigcap_{j\not\in I\cup
J}U^c_j
\bigcap_{i\in I-J} (U_i-U_i^\circ);\\
X_{I,I\cup J}&=& \bigcap_{i\in I\cap J}U_i\bigcap_{j\not\in
I}U^c_j
\bigcap_{k\in I-J} (U_k-U_k^\circ);\\
X_{J,I\cup J}&=& \bigcap_{i\in I\cap J}U_i\bigcap_{j\not\in
J}U^c_j \bigcap_{k\in J-I} (U_k-U_k^\circ),
\end{eqnarray*}
we have
$$
X_{I\cup J, I}\cap X_{I\cup J,J}= X_{I\cap J,I}\cap X_{I\cap J,J}=
\bigcap_{i\in I\cap J}U_i \bigcap_{j\not\in I\cup J}U_j^c
\bigcap_{k\in I\cup J-I\cap J} (U_k-U_k^\circ).
$$
This says, by \eqref{eqn_2.3.2},
$$
X_{I\cup J, I}\cap X_{I\cup J,J}= X_{I\cap J,I}\cap X_{I\cap J,J}=
X_{I\cup J,I\cap J}.
$$
It is also easy to check that they are same as $X_I\cap X_J$. This
coincides with the definition of $X_{I,J}=X_{J,I}$. These imply
(P1)-(P5) in Definition \ref{defn_2.2.1}.
 q.e.d.
\subsection{Language of germs}\label{sect_2.4}
Let $(\mc X,\Phi)$ be a virtual manifold and $\mathbf X$ be its
virtual space. $\mathbf X$ admits a partition. For each $I\in \mc
N$ we define
$$
\mathbf X_I= [X_I]-\bigcup_{J\subset I}[X_J].
$$
Then clearly $\{\mathbf{X}_I\}$ is a partition of $\mathbf X$.
Note that $[x]\in \mathbf{X}_I$ if and only if the support of
$[x]$ is $X_I$. For $[x]\in \mathbf{X}_I$, let $x$ be the
corresponding point in $X_I$. Let $\mc N_{[x]} $ be the set of
open neighborhoods of $x$ in $X_I$ given as following:
$$
\mc N_{[x]}=\{\phi_I\inv(V)| V \mbox{ is neighborhood of } [x]\}.
$$
Suppose that  $U$ is an element in $\mc N_{[x]}$.  Then for any
$y\in U$ we have a unique element in $\mc N_{[y]}$: if the support
of $[y]$ is $X_I$, we take the element to be $U$; otherwise,
suppose the support of $[y]$ is $X_J, J\subset I$, we take the
element to be $\phi_{I,J}(U\cap X_{I,J})$.

We wish to develop the theory on $\mathbf X$ via structure $(\mc
X,\Phi)$. Let
$$
\Gamma=\{U|U \in \mc N_{[x]} \mbox{ for some } [x]\}.
$$
We say such a collection $\Gamma$ is a {\em complete collection } if
\begin{itemize}
\item $\{[U]|U\in \Gamma\}$ covers $\mathbf X$; and
\item for any $U=U_{[x]}\in \mc N_{[x]}$, the induced element $U_{[y]}$
for $[y]\in [U]$ is also in $\Gamma$.
\end{itemize}
Given a complete collection $\Gamma$, we set
$$
X'_I=\bigcup_{U\subset X_I, U\in \gamma} U.
$$
Then $X'_I\subseteq X_I$. Set
$$
X'_{I,J}=X'_I\cap X_{I,J},
$$
and $\phi'_{I,J}=\phi_{I,J}|_{X'_{I,J}}$. It is easy to see that
$(\mc X',\Phi')$ forms a virtual manifold.
 Moreover $\mathbf
X'=\mathbf X$. Let $\mc Z(\mc X,\Phi)$ be the collection of virtual manifolds
constructed by this way. $\mc Z(\mc X,\Phi)$ admits a partial order: let $(\mc
X',\Phi')$ and $(\mc X'',\Phi'')$ be two virtual manifolds in $\mc
Z(\mc X,\Phi)$, we say that $(\mc X',\Phi')\prec (\mc X'',\Phi'')$ if and only if
for any $I$
$$
X'_I\subseteq X''_I.
$$
Note that for any two virtual manifolds $(\mc X_i,\Phi_i)\in \mc
 Z(\mc X,\Phi),i=1,2,$ there exists an  $(\mc X_3,\Phi_3)\in
 \mc Z(\mc X, \Phi) $ such that
 $$
(\mc X_3,\Phi_3)\prec (\mc X_i,\Phi_i), i=1,2.
 $$
 The germ of $(\mc X,\Phi)$, denoted as $(\mc
X,\Phi)^{germ},$ is defined to be the direct limit of $\mc Z$. We
propose the principle of theory on virtual manifolds: \v {\em A
theory on $(\mc X,\Phi)$ is a theory on $(\mc X,\Phi)^{germ}$.} \v
By the spirit of germs, a theory $P$ on $(\mc X,\Phi)^{germ}$ is
constructed on some $(\mc X',\Phi')\in \mc Z(\mc X,\Phi)$.

\begin{remark}\label{rmk_2.4.1}
Sometimes, we need a more subtle version. A theory $P$ on $(\mc
X,\Phi)^{germ}$ means
\begin{itemize}
\item theory $P_{[x]}$ on some $U_{[x]}$ for each $[x]$;
\item
theory $P_{[x]}$ is compatible with theory $P_{[y]}$ for any $[x]$
and $[y]$.
\end{itemize}
The compatibility is equivalent to the following statement: let
$P_{[x]}$ be the theory on $U_{[x]}$; for any $y\in U_{[x]}$,
$U_{[x]}$ induces an element $U'\in \mc N_{[y]}$; then we require
that the theory $P_{[y]}$  and $P_{[x]}$ are compatible  on
$U'\cap U_{[y]}$.
\end{remark}
We remark that for any collection
$$
\Gamma'=\{U_{[x]}|U_{[x]}\in \mc N_{[x]}\}_{[x]\in \mathbf X}
$$
we can generate a complete collection $\Gamma$ that contains $\Gamma'$.


\subsection{Category of virtual manifolds}\label{sect_2.5}

We construct the category of virtual manifolds. The main task is
to construct the maps between virtual manifolds.

Let $(\mc X,{\Phi})$ and $(\mc B,\Psi)$ be two virtual manifolds,
 $\mathbf X$ and $\mathbf B$ be their virtual spaces. Let
$\mathbf{f}: \mathbf X\to \mathbf B$ be a continuous map. We
define lifts of $\mathbf f$ on $(\mc X, \Phi)^{germ}\to (\mc
B,\Psi)^{germ}$ in terms of language of germs in the sense of
Remark \ref{rmk_2.4.1}.
\begin{defn}\label{defn_2.5.1}
A collection of maps $\mc F=\{f_{[x]}:U_{[x]}\to V_{ \mathbf
f([x])}\},$ where $U_{[x]}\in \mc N_{[x]}$ and $V_{f_{[x]}}\in \mc
N_{\mathbf ([x])}$, is a lift of $\mathbf f$  if $f_{[x]}$ is
compatible with $f_{[y]}$ for all $[x]$ and $[y]$, namely,
\begin{enumerate}
\item[]  {\em Case 1, }if $U_{[x]}$ and $U_{[y]}$ are in same
$X_I$, $f_{[x]}=f_{[y]}$ on $U_{[x]}\cap U_{[y]}$; \item[] {\em
Case 2,} otherwise, suppose $U_{[y]}\subset X_J, J\subset I$, then
$f_{[x]}$ is a lifting of $f_{[y]}$ on the fibration
$$
\phi_{I,J}\inv(\phi_{I,J}(U_{[x]}\cap X_{I,J})\cap U_{[y]})\to
\phi_{I,J}(U_{[x]}\cap X_{I,J})\cap U_{[y]}.
$$
\end{enumerate}

Two lifts $\mc F'$ and $\mc F''$ are equivalent if for all $[x]$,
$f'_{[x]}$ and $f''_{[x]}$ are compatible. An equivalence class
$[\mc F]$ is called a lift of $\mathbf f$ on $(\mc
X,\Phi)^{germ}$. $(\mathbf f,[\mc F])$ is called a {\em virtual
map} from $(\mc X,\Phi)$ to $(\mc B,\Psi)$.
\end{defn}
Note that a collection of maps between $(\mc X,\Phi)$ and $(\mc
B,\Psi)$ given by
$$
f_I: X_I\to B_I
$$
satisfying that $f_J$ is a lifting of $f_I$ on the fibration
$$
\Phi_{J,I}: X_{J,I}\to X_{I,J}, I\subset J
$$
determines a virtual map.
\begin{prop}\label{prop_2.5.1}
Let
\begin{eqnarray*}
&&(\mathbf{f}_1,[\mc F_1]): (\mc X_0,\Phi_0)\to
(\mc X_1,\Phi_1),\\
&&(\mathbf{f}_2,[\mc F_2]): (\mc X_1,\Phi_1)\to (\mc X_2,\Phi_2);
\end{eqnarray*}
be two virtual maps, then they are composed to a (unique) virtual
map
$$
(\mathbf{f}_2\circ\mathbf{f}_1,[\mc G]): (\mc X_0,\Phi_0)\to (\mc
X_2,\Phi_2).
$$
\end{prop}
{\bf Proof. } This is obvious via the definition of germs. q.e.d.
\v
\begin{defn}\label{defn_2.5.2}
$(\mathbf{f},[\mc F])$ is called a smooth map if all lifts
$f_{[x]}$ are smooth.
\end{defn}
\n Then, by Proposition \ref{prop_2.5.1}, we make the following
definition.
\begin{defn}\label{defn_2.5.3}
The category of (smooth) virtual manifolds is denoted by
${\mc{V}}$. The objects of the category are virtual manifolds. The
morphisms between two virtual manifolds are equivalence classes of
smooth virtual maps.
\end{defn}

\section{Integration theory on virtual manifolds}\label{sect_3}

We  develop the integration theories on virtual manifolds. A
similar theory  holds for virtual orbifolds.

\subsection{Forms on virtual manifolds}\label{sect_3.1}

Let $(\mc X,\Phi)$ be a virtual manifold.
\begin{defn}\label{defn_3.1.1}
A pre-$k$-form on $(\mc X,\Phi)$ is
$$
\alpha=\{\alpha_I\in \Omega^k(X_I)|I\in \mc N\}
$$
such that
$$
\alpha_J=\phi^\ast_{J,I}\alpha_I
$$
on $X_{J,I}$. Two pre-$k$-forms $\alpha'$ and $\alpha''$,  on
$(\mc X',\Phi')$ and $(\mc X'',\Phi'')$ respectively, are
equivalent if they admits on a smaller $(\mc X''',\Phi''')\in \mc
Z(\mc X,\Phi)$. Let $[\alpha]$ be the equivalence class. It is called a
$k$-form on $(\mc X,\Phi)$.
\end{defn}
Namely, a form on $(\mc X,\Phi)$ is a form on $(\mc
X,\Phi)^{germ}$.
 Let $\Omega^k(\mc X)$ be set of $k$-forms on $(\mc
X,\Phi)$. Then $(\Omega^\ast(\mc X),d)$ is a complex. Define
$$
H^\ast_{dR}(\mc X)= H^\ast (\Omega^\ast(\mc X),d).
$$

We next consider a very different type of forms on virtual
manifolds. Let $\Theta_{J,I}$ be the Thom forms of the bundle
$\Psi_{J,I}: X_{J,I}\to X_{I,J}$. To avoid the unnecessary
complication caused by the degree of forms, we always assume that
the degree of $\Theta_{J,I}$ is even.
\begin{defn}\label{defn_3.1.2}
A set of  forms $\Theta= \{\Theta_{J,I}\}_{I\subseteq J} $ is
called a {\em transition data} of $\mc X$ if it satisfies the
following compatibilities: for any $I$ and $J$,
$$
\Theta_{I\cup J, I\cap J} = \Psi^*_{I\cup J,I}\Theta_{I,I\cap
J}\wedge \Psi^*_{I\cup J,J} \Theta_{J,I\cap J}
$$
on $X_{I\cup J,I\cap J}$.
\end{defn}
\begin{defn}\label{defn_3.1.3}
A pre-virtual form on $(\mc X,\Phi)$ is
$$
\mfk z=\{z_I\in \Omega^\ast(X_I)|I\in \mc N\}
$$
such that
$$
z_J=\phi^\ast_{J,I}z_I\wedge \Theta_{J,I}
$$
on $X_{J,I}$ for some transition data $\Theta$. $\mfk z$ is called
a $\Theta$-form on $\mc X$. Two pre-$k$-forms $\mfk z'$ and $\mfk
z''$,  on $(\mc X',\Phi')$ and $(\mc X'',\Phi'')$ respectively,
are equivalent if they admits on a smaller $(\mc X''',\Phi''')\in
\mc Z$. Let $[\mfk z]$ be the equivalence class. It is called a
virtual form on $(\mc X,\Phi)$. The virtual degree of $\mfk z$ is
the degree of $z_\emptyset$.
\end{defn}
Let $\Omega^\ast_v(\mc X)$ be set of virtual forms on $(\mc
X,\Phi)$. Then $(\Omega^\ast_v(\mc X),d)$ is a complex. Define
$$
H^\ast_{v,dR}(\mc X)= H^\ast (\Omega^\ast_v(\mc X),d).
$$

We may define the support of a form or a virtual-form. Let us take
a form $[\alpha]$ as an example. Suppose $\alpha= (\alpha_I)$. For
any $[x]\in \mathbf X$, let $X_I$ be its support, we say
$$
[x]\in \mathrm{supp}([\alpha]) \iff x\in \mathrm{supp}(\alpha_I).
$$
If $[\alpha]$ is compact supported in $\mathbf X$ (or $\mathbf
X^\circ$), we write $[\alpha]\in \Omega^*_c(\mc{X})$ (or
$\alpha\in \Omega^*_c(\mc X^\circ)$). Similarly, we can define
$\Omega_{v,c}(\mc X)$ and $\Omega_{v,c}(\mc X^\circ)$. For most of
time, we are interested in forms  $\Omega_{v,c}(\mc X^\circ)$.

If $\alpha\in \Omega^*(\mc X)$ and $\mfk z\in \Omega_{v,c} (\mc
X)$, $\alpha\wedge \mfk z$ is in $\Omega_{v,c}(\mc X)$.

\subsection{Integrations on virtual manifolds.}\label{sect_3.2}
We now describe how to define $ \int_{\mathbf X} [\mfk z] $ for
$[\mfk z]\in \Omega_{v,c}(\mc X)$, where $\deg\mfk z=\dim \mc X $.
Suppose $\mfk z=(z_I)$ is a $\Theta$-form on $(\mc X,\Phi)$
representing $[\mfk z]$. The definition is almost obvious because
of the following reason: suppose $[V]\subseteq [X_I]\cap[X_J]$,
let $V_I=\pi^{-1}_I([V]), V_J=\pi^{-1}_{J}([V]), V_{I\cap J}=
\pi^{-1}_{I\cap J}([V])$ and $V_{I\cup J}=\pi^{-1}_{I\cup
J}([V])$, then
$$
\int_{V_I} z_I= \int_{V_{I\cup J}}z_{I\cup J} =\int_{V_J} z_J.
$$
The equalities of two ends are due to the Thom isomorphism.
Similarly, the middle term can also be replaced by $\int_{V_{I\cap
J}}z_{I\cap J}$. Hence, $\int_X\mfk z$ is well defined on
$[X_I]\cap [X_J]$. Now define
\begin{equation}\label{eqn_3.2.1}
\int_{\mathbf X} \mfk z= \sum_I\int_{[X_I]}z
-\sum_{I,J}\int_{[X_I]\cap [X_J]} z
+\sum_{I,J,K}\int_{[X_I]\cap[X_J]\cap[X_K]} z-\cdots.
\end{equation}
It is not hard to see that $\int_{\mathbf X}[\mfk z]$ is
independent of choice of representatives of $[\mfk z]$.

Let
$$
\Gamma=\{ U| U\in \mc N_{[x]} \mbox{ for some }[x]\}
$$
be a collection of sets such that $\{[U]|U\in \Gamma\}$ covers
$\mathbf X$. Let
$$
X'_I=\cup_{U\subset X_I} U.
$$
Then
\begin{equation}\label{eqn_3.2.1}
\int_{\mathbf X} \mfk z= \sum_I\int_{[X_I']}z
-\sum_{I,J}\int_{[X_I']\cap [X_J']} z
+\sum_{I,J,K}\int_{[X_I']\cap[X_J']\cap[X_K']} z-\cdots.
\end{equation}

Furthermore, we can construct a new virtual manifold out of
$\{X'_I\}$. Recall that $\{[X'_I]\}$ forms a cover of $\mathbf X$.
We now apply Example \ref{ex_2.3.1}. To do this, we re-index
the index set
$$\iota:\mc N\to N^{[1]}=\{0,1,\ldots, 2^N-1\}
$$
by requiring $\iota(\emptyset)=0$. Rewrite the sets $\{[X_I']\}$
as $\{Y_i\}$, i.e, $Y_{\iota(I)}= [X'_I]$. Set $\mc N^{[1]}
=2^{N^{[1]}}$. Then as explained in Example \ref{ex_2.3.1},
we can construct a patchable pair
$(\mc Y, \mathbf \Psi)$. We now construct a virtual manifold
from this pair. For $I^{[1]}=\{i_1,\ldots, i_k\}\in \mc N^{[1]}$,
we have a set
$$
\mathbb{I}=\{ I_{i_1},\ldots, I_{i_k}\}, \mbox{ where }
I_{i_j}=\iota\inv(i_j).
$$
Set
$$
I_{\max}= \bigcup_{j=1}^k I_{i_j}\in \mc N.
$$
We define
$$
Z_{I^{[1]}}= \pi_{I_{\max}}\inv(Y_{I^{[1]}}).
$$
Note that $Y_{I^{[1]}}\subset \mathbf X$,  in particular, in
$[X_{I_{\max}}]$. Hence $Z_{I^{[1]}}$ is in $X_{I_{\max}}$.
For $I^{[1]}\subset J^{[1]}$,
$$
Z_{I^{[1]},J^{[1]}}= \pi_{I_{\max}}\inv(Y_{I^{[1]},J^{[1]}}),
\;\;
Z_{J^{[1]},I^{[1]}}= \pi_{J_{\max}}\inv(Y_{J^{[1]},I^{[1]}}).
$$
and
$$
\Psi_{J^{[1]},I^{[1]}}:Z_{J^{[1]},I^{[1]}}\to Z_{I^{[1]},J^{[1]}}
$$
is induced from $\Phi_{J_{\max},I_{\max}}$.
Set
$$
\mc Z_\Gamma=\{Z_{I^{[1]}}\}, \;\;\;
\mathbf\Psi=\{\Psi_{J^{[1]},I^{[1]}}\}.
$$
It is straightforward
to prove that
\begin{prop}\label{prop_3.2.1}
$(\mc Z_\Gamma,\mathbf\Psi)$ is a virtual manifold and its
virtual space $\mathbf Z_\Gamma$ is same as $\mathbf X$.
\end{prop}
Note that the integration on $\mathbf Z_\Gamma$ is same
as on $\mathbf X$. We call this new virtual manifold
$\mc Z_\Gamma$ to be a modification of $\mc X$.

To manipulate integrations, it is convenient to develop some type of
 theorems of
partition of unity on virtual manifolds. These are discussed in the next
sub-section.

\subsection{Partition of unity}\label{sect_3.3}

Let $(\mc
X,\Phi)$ be a virtual manifold. For simplicity, we assume that
$\partial \mc X=\emptyset$. Also we assume \eqref{eqn_3.2.2}
holds.

\begin{defn}\label{defn_3.3.1}
Let $\mc W=(W_I)$, where $W_I$'s are open subset of $X_I$. We say
$\mc W$ is {\em pre-compact} if it satisfies:
\begin{itemize}
\item[a1.] $\overline{[{\mc W}]}$ is compact in $X$, where
$$
[\mc W]:=\cup [W_I] ;
$$
\item[a2.] for each $x\in [X_I]$, $\overline{\pi_I\inv\cap W_I}$
is compact; \item[a3.] $\Psi_{J,I}(W_{J,I})= W_{I,J}$, where
$W_{J,I}=X_{J,I}\cap W_J$.
\end{itemize}
\end{defn}
Let $\mc W=(W_I)$ be pre-compact. Then we have
\begin{lemma}\label{lemma_3.3.1}
There exists open subset $Y_I\subset [X_I]$ for each $I$ such that
\begin{itemize}
\item[b1.] $\overline{[\mc W] }\subset \cup Y_I$; \item[b2.] $\bar
Y_I\subset [X_I]$.
\end{itemize}
\end{lemma}
{\bf Proof. } For any $[x]\in \overline{[\mc W] }$, let $X_K$ be
the support of $[x]$, i.e, there exists a unique $x\in X_K$ such
that $\pi_K(x)=[x]$. Hence there exists a small neighborhood
$B_K(x) \subset X_K$ such that $\bar B_K(x)\subset X_K$. Then
$\pi_K(\bar B_K(x))\subset [X_K]$.
$$
\{[B_K(x)],x\in \overline{[\mc W] }\}
$$
covers $\overline{[\mc W] }$. Since the latter one is compact,
there exists a finite cover, denoted by $\{[B_{K_i}(x_i)], 1\leq
i\leq l\}$ for some $l<\infty$. Set
$$
Y_I= \bigcup_{K_i=I} [B_{K_i}(x_i)].
$$
Clearly, the lemma follows. q.e.d.

\v Set $Z_I=\pi_I\inv Y_I$, $Z_I$ is an open subset of $X_I$.

\begin{defn}\label{defn_3.3.2}
We say $\{\eta_I\}$ is a  smooth partition of unity with respect to
$\mc W$ (i.e, for any $x\in [\mc W], \sum \eta_I(x)=1$) if
\begin{itemize}
\item[c1.] $\mathrm{supp}(\eta_I)\subset Y_I$; \item[c2.]
$\beta_I=\pi_I^\ast\eta_I$ is smooth on $Z_I$.
\end{itemize}
\end{defn}
It is not clear that   a {\em smooth}
 partition of unity exists on a virtual manifold. But we can prove the
 existence on its modification. Let $\mc Z_\Gamma$ be a modification
 of $\mc X$ explained in the end of last subsection. Correspondingly,
 we have
 $\mc W$ on $\mc Z_\Gamma$.

\begin{lemma}\label{lemma_3.3.2}
For a proper choice of $\Gamma$, there exists a partition of unity $\{\eta_{I^{[1]}}\}$
on $\mc Z_\Gamma$ with respect to
$\mc W_\Gamma$.
\end{lemma}
{\bf Proof. } Step 1, by classic result, there is a continuous partition unity
${\eta_I}$ on $\mc X$ with respect to $\mc W$. We can first construct
it on topological space $\mathbf X$, then pull back functions to
$X_I$.
Set
$$
[S_I]=\bigcup_{J\prec I}\partial\overline{[X_J]}
$$
and
$S_I=\pi_I\inv([S_I])$.
It is not hard to construct   $\eta_I$ such that they are smooth away
from $S_I$.
\v
Step 2, by modifying  functions $\eta_I$, we may have
$\eta_I'$ defining on $X'_I=\pi_I\inv A_I$, where
$$
A'_I\subset \bar A'_I\subset  [X'_I]
$$
and $\{A'_I\}$ covers $\mathbf X$.

This can be done inductively. We start with $I=\emptyset$. Set
$\eta_I'=\eta_I$. Next let $I=\{i\}, 1\leq i\leq N$. We mollify
$\eta_I$ at $S_I$ and get a smooth function $\eta_I'$ on $X'_I$.
We can then modify $\eta_J$ on $X_J$ such that
$\eta_J$ equals $\eta_I$ on the overlapping area for  $I\prec J$.
Inductively, we do the same construction for $I=\{i,j\}$ and so on.

We conclude that for $I\prec J$, $\eta_I'$ matches $\eta_J'$ on the overlapping area away from
$S_J$. Now we can choose proper $A'_I$ such that
$\eta_I'$ matches $\eta_J'$ on the $A'_I\cap A'_J$. Moreover we can assume
that $\{A'_I\}$ covers $\mathbf X$.
\v
Step 3, Set $\Gamma=\{ A'_I\}$. Then we can construct a
virtual manifold $\mc Z_\Gamma$ with smooth functions $\eta'_{I^{[1]}}$
on $Z_{I^{[1]}}$. Now on $\mathbf Z_\Gamma$, set
$$
\eta_{I^{[1]}}=\frac{\eta_{I^{[1]}}'}{\sum_{J^{[1]}} \eta_{J^{[1]}}' }.
$$
This gives a smooth partition of unity on $\mc Z_\Gamma$. q.e.d.

\v\n
Suppose $\mc X$ has a partition of unity.
 Another way to get $\int_X \mfk z$ is to use partition of
unity. Set
$$
W_I=(\mathrm{supp} (z_I))^\circ
$$
and $\mc W=\{W_I\}$. Let $\{\eta_I\}$ be a partition of unity with
respect to $\mc W$. Then equivalently, we define the integral
\eqref{eqn_3.2.1} to be
\begin{equation}\label{eqn_3.3.2}
\int_X \mfk z= \sum_I\int_{X_I}\eta_I z_I.
\end{equation}

\subsection{The Stokes' theorem on virtual manifolds}\label{sect_3.4}
In order to define invariants, we  focus on the following case:
(1), $\mfk z\in \Omega_{\Theta,c}^{k}(\mc{X^\circ})$ is close;
(2), $a\in \Omega^{d-k}(\mc{X})$ is close. Define
$$
\mu_{\mfk z}(a)= \int_{X} a\wedge\mfk z.
$$
Then,  Stokes' theorem and the Thom isomorphism imply that
\begin{lemma}\label{lemma_3.4.1}
Suppose $a$ and $b$ are exact,
$$
\mu_{\mfk z}(a)=\mu_{\mfk z}(b).
$$
\end{lemma}
{\bf Proof: } Suppose $dc= a-b$, where $c=(c_I)$. Then
$$
\mu_{\mfk z}(a)- \mu_{\mfk z}(b) =\sum_I\int_{  X_I} dc_I\wedge
\eta_Iz_I =-\sum_I \int_{X_I}d(\eta_I)c_I\wedge z_I= 0.
$$
The last equality follows from the Stokes' theorem, the Thom
isomorphism, and that $\{\eta_I\}$ is a partition of unity. q.e.d.
\v \n A similar argument implies  the Stokes' theorem for virtual
manifolds. Suppose $\mfk z\in \Omega_{\Theta,c}(\mc L)$. The
restriction of $\mfk z$ on $\partial \mc L$, denoted by
$i^\ast\mfk z$, is a form in $\Omega_{\Theta,c}(\partial \mc L)$.
Here $i:\partial \mc L\to \mc L$ is the standard embedding. Then
\begin{theorem}[Stokes' Theorem]\label{thm_3.4.1} For $\mfk z\in \Omega_{\Theta,c}(\mc L)$
$$
\int_X d\mfk z= \int_{\partial X}i^\ast\mfk z.
$$
\end{theorem}
As a consequence,
\begin{corollary}\label{cor_3.4.1}
For any close form $\mfk z\in \Omega_{\Theta,c}(\mc L)$
$$
\int_{\partial X}i^\ast\mfk z=0.
$$
\end{corollary}

Note that we just have a pairing
$$
\mu: H^\ast_{v,c}(\mc X^\circ)\times H^\ast(\mc X)\to \mathbb{R}.
$$

\subsection{Virtual bundles and the Euler classes}\label{sect_3.5}

Let $(\mc X,\Phi)$ be a virtual manifold without boundary.
\begin{defn}\label{3.5.1}
A virtual bundle  over $(\mc X,\Phi)$ is $\mc E=(E_I)$, where
$E_I\to X_I$ is a vector bundles, such that for $I\subset J$
$$
E_J|_{X_{J,I}}= \phi_{J,I}^\ast (X_{J,I}\oplus E_I|_{X_{I,J}}).
$$

A section of the bundle is $\mc S=(S_I)$, where
$S_I: X_I\to E_I$ is a section of bundle, such that for $I\subset J$
$$
S_J(x,v)= (v,S_I(x)).
$$
Here $(x,v)\in X_{J,I}$ is the local coordinate. The section is transverse to
0-section if each $S_I$ is transverse.

Let $\Theta=(\Theta_{J,I})$ be a transition data.
A $\Theta$-Thom form of $\mc E$ is $\Lambda=(\Lambda_I)$, where $\Lambda_I$ is a Thom form
of $E_I$, such that for $I\subset J$
$$
\Lambda_J= \Theta_{J,I}\wedge \Lambda_I
$$
on $E_J|_{X_{J,I}}$. $S^\ast \Lambda$ defines a pre-$\Theta$-form. $[S^\ast \Lambda]$
is called an Euler form of $\mc E$.
\end{defn}

The proofs of the  following two statements are straightforward. We leave the
proofs to readers.
\begin{prop}\label{prop_3.5.1}
Let $\mc E=(E_I)$ be a virtual bundle over a virtual manifold
$(\mc X,\Phi)$. Then there exists a sub-virtual manifold $(\mc
X',\Phi')$ in the sense of Remark \ref{rmk_3.2.1} such that, over
the virtual bundle given by $\mc E'=(E_I|_{X_I'})$ we have
\begin{itemize}
\item[1,] a section $S$,
\item[2,] a $\Theta$-Thom form $\Lambda$, and therefore
\item[3,] an Euler class $S^\ast \Lambda$.
\end{itemize}
\end{prop}
Let $\mc E\to \mc X$ be a virtual bundle with a virtual transverse section $S$. Then
$S\inv(0)$ is a virtual manifold.
\begin{prop}\label{prop_3.5.2}
Let $\mc E^1$ and $\mc E^2$ be  two virtual bundles over $(\mc
X,\Phi)$. Let $\Lambda^i$ be their $\Theta$-Thom forms. Let
$S^i,i=1,2,$ be two transverse sections of $\mc E^i$. Then
$S=S^1\oplus S^2$ is a transverse section of $\mc E^1\oplus \mc
E^2$. Set $\lambda=\Lambda^1\wedge\Lambda^2$. We have
$$
\int_{\mc X} a\wedge (S^\ast \Lambda)
=\int_{(S^1)\inv(0)} a\wedge (S^2)^\ast(\Lambda^2)
=\int_{(S^2)\inv(0)} a\wedge (S^1)^\ast(\Lambda^1).
$$
\end{prop}

\section{$G$-virtual manifolds  and  localization }\label{sect_4}

\subsection{$G$-virtual manifolds, equivariant forms and integration}\label{sect_4.1}
The discussion given in the previous section can be generalized to
the equivariant case. Let $G$ be a compact Lie group.
\begin{defn}\label{defn_4.1.1}
By a {\em $G$-virtual manifold} ($\mc{X},\Phi)$,  we mean that
(a.) $(\mc X,\Phi)$ is a virtual manifold, (b.) each $X_I$ is
$G$-manifold and (c.) $\Psi_{J,I}: X_{J,I}\to X_{I,J}$ are
$G$-equivariant bundles for any $I\subset J$.
\end{defn}
To study the $G$-equivariant integration theory on $\mc X$, we
may consider  {\it $G$-equivariant transition data}
$\Theta_G=\{\Theta^G_{J,I}\}_{I\subseteq J}$. Then similarly, we
may define:
 $G$-equivariant
forms $\Omega_{G}^*(\mc X)$, $G$-equivariant $\Theta_G$ forms
$\Omega_{\Theta_G}(\mc X)$, $\Omega_{\Theta_G,c}(\mc X)$ and
$G$-invariant partition of unity, etc. For $\zeta=(\zeta_I)\in
\Omega_{\Theta_G,c}(\mc X)$, as before, we define
$$
\int^G_X \zeta= \sum_I\int_{ X_I} \eta_I\zeta_I.
$$
For a $G$-closed form $\zeta\in \Omega_{\Theta_G,c}(\mc X^\circ)$
and $\alpha\in \Omega^*_G(\mc X)$, define
$$
\mu_{\zeta}(\alpha)= \int_X^G \alpha\wedge \zeta.
$$

For any space $Y$ with $G$ action, the notation $Y^G$ denotes the
fix loci of the action. Let
$\mc X^G=\{X_I^G\}$
Then
\begin{lemma}\label{lemma_4.1.1}
$\mc X^G$ is a virtual (sub-)orbifold (of $\mc X$).
 \end{lemma}
\n This follows directly from the definition. We skip the
 proof.
\v We will discuss the  abelian localization formula of
Atiyah-Bott type for the integration $\mu_{\zeta}(\alpha)$. For
simplicity, we assume that $G=S^1$. We also assume that $\mc X$
is compact, boundary free for the sake of Stokes' theorem.
Otherwise,  the compact-supportedness of $\zeta$ would take care
the issue.

\subsection{Virtual normal bundles and their Euler classes}
\label{sect_4.2}

For each $X_I^G$, let $N_I$ denote its $G$-normal bundle in $X_I$.
Now focus on $X_{J,I}\to X_{I,J}, I \subseteq J$. By restricting
on $X_{I,J}^G$, this bundle splits as
\begin{equation}\label{eqn_4.2.1}
X_{J,I}|_{X_{I,J}^G}=X_{J,I}^G\oplus P
\end{equation}
for some $G$-invariant bundle $P$. Therefore,
$$
N_J|_{X_{J,I}^G}= \Psi_{J,I}^\ast (N_I\oplus P).
$$
This says that
\begin{lemma}\label{lemma_4.2.1}
$\mc N=\{N_J\}$ forms a virtual bundle over $\mc X^G$.
\end{lemma}
Let $e_{J,G}$ be the equivariant Euler forms of $N_J$ over
$X_J^G$. We may arrange them such that
\begin{equation}\label{eqn_4.2.2}
e_{J,G}|_{X_{J,I}^G} =\Psi_{J,I}^\ast e_{I,G}\wedge \Theta_G(P),
\end{equation}
where $\Theta_G(P)$ is the equivariant Thom form on $P$.
 Denote $\{e_{I,G}(X_I^G)\}$
by $e_G(X^G)$.

On the other hand, Let $\Theta_G=(\Theta_{I,J}^G)$ be the
transition data. According to \eqref{eqn_4.2.2}, we may assume
that
\begin{equation}\label{eqn_4.2.3}
\Theta_{J,I}^G= \Theta_G(X_{J,I}^G) \wedge \Theta_G(P).
\end{equation}
over $X_{J,I}|_{X_{I,J}^G}$. This induces a transition data
$$
\tilde\Theta_G=\{\Theta_G(X_{J,I}^G)\}
$$
on $\mc X^G$.

For any $\Theta_G$-form $\zeta=(\zeta_I)$, we find that
$$
\frac{\zeta_J}{e_{J,G}}
$$
forms a $\tilde\Theta$-form on $X^G$:
$$
\frac{\zeta_J}{e_J^G} =\Psi_{J,I}^\ast
\frac{\zeta_I\wedge\Theta_{J,I}^G}{e_{I,G} \wedge \Theta_G(P)}
=\Psi_{J,I}^\ast
\frac{\zeta_I\wedge\tilde\Theta_G(X_{J,I}^G)}{e_{I,G}}
$$
Hence, we conclude that
\begin{prop}\label{prop_4.2.1}
For a $\Theta^G$ form $\zeta$, $\{\zeta_I/e_{I,G}\}$ forms a
$\tilde\Theta$-form on $X^G$. The form is denoted by
$e_\zeta(X^G)$.
\end{prop}

\subsection{The Abelian localization formula}\label{sect_4.3}
The standard localization technique  implies that
\begin{theorem}\label{thm_4.3.1}
Let $\mc X$ be a finite dimensional virtual manifold with $G=S^1$
action. Let $X$ be its virtual space. Let $\zeta\in
\Omega_{\Theta_G,c} (\mc X^\circ)$ and $\alpha\in
\Omega^*_G(\mc X)$, then
$$
\mu_{\zeta}(\alpha)=\int_{X^G}\frac{ i_{X^G}^*(\alpha\wedge
\zeta)}{ e_{G}(X^G)}.
$$
\end{theorem}
The right hand side in the formula can be thought as an
integration on virtual manifold $\mc X^G$:
$$
\mu_{e_\zeta(X^G)}(i_{X^G}^\ast \alpha)
$$
\n {\bf Proof: } Let $\Omega_I,\Omega_I'$ be two equivariant Thom
forms on $N_I$. By identifying $N_I$ with a  neighborhood $\tilde
U(X^G_I)$ of $X_I^G$, we require that they are equal in a smaller
neighborhood $U(X_I^G)\subseteq \tilde U(X^G_I)$ and
$i_{X_I}^*\Omega_I= e_{I,G}(X_I^G)$. Moreover, we require that:
(1) the support of $\Omega_I'$ is contained in that of $\Omega_I$;
(2) on the overlapping area $ N_J\cap X_{J,I} \to N_I\cap X_{I,J},
$
$$
\Omega_J=\Omega_I\wedge \Theta^G_{J,I}, \mbox{ and }
\Omega'_J=\Omega'_I\wedge \Theta^G_{J,I}.
$$
It is not hard to see that such pairs always exist.

Then by the Thom isomorphism,
$$
\sum_I \int_{X_I}\frac{\eta_I\alpha_I\wedge \zeta_I\wedge
\Omega_I'}{\Omega_I} =\sum_I \int_{ X_I^G}
\frac{i^*_{X_I^G}(\eta_I\alpha_I\wedge \zeta_I)}{e_G(X_I^G)},
$$
Note that the right hand side is same as the right hand side of
the formula in the theorem.

 On
the other hand, let
$$
\tilde{\alpha}_I= \alpha_I- \alpha_I\wedge
\frac{\Omega_I'}{\Omega_I}.
$$
Then $\tilde\alpha=(\tilde \alpha_I)$ becomes a $G$-equivariant
form  supported away from $X^G$. It remains to prove that
$$
\mu_{\zeta}(\tilde\alpha)=0.
$$
The proof is standard. In fact, there exists a form
$\beta=(\beta_I)$ supported away from $X^G$ such that
$\tilde\alpha_I= d_G \beta_I$. Then
$$
\mu_\zeta(\tilde\alpha) =\sum_I\int_{X_I}\eta_I d\beta_I
=-\sum_I\int_{X_I} d(\eta_I) \beta_I=0
$$
For the last equality, we use the fact that $\sum_I\eta_I=1$.
q.e.d.

\section{Fredholm systems and Stabilizations}\label{sect_5}

\subsection{The Fredholm set-up}\label{sect_5.1}
We start with the following set-up.
\begin{defn}\label{defn_5.1.1}
A Fredholm system consists of following data:
\begin{enumerate}
\item[(B1)] let $\pi:\mcf\to \mcb$ be a Banach orbifold bundle
over a Banach orbifold $\mcb$; \item[(B2)] let $S: \mcb\to \mcf$
be a proper smooth section.
 In particular, the
properness implies that $M=S^{-1}(0)$ is compact; \item[(B3)] for
any $x\in M$, let $L_x$ be the linearlization of $S$ at $x$
$$
L_x: T_x\mcb\to \mcf_x.
$$
We assume that $L_x$ is a  Fredholm operator. Let $d$ be the index
of the operator.
\end{enumerate}
We refer the triple $(\mcb,\mcf,S)$ as a {\em Fredholm system}.
$M$ is called  the {\em moduli space} of the system.
\end{defn}

A core topic in studying moduli problems is to define  invariants
on such a system. This is  based on the study of $M$. It is well
known that if $L_x$ is surjective for all $x\in M$, $M$ is a
compact smooth orbifold. Then $M$ can be thought as a cycle in
$H_d(\mcb)$ representing the Euler class of  bundle $\mcf\to\mcb$.
Let $a\in H^d(\mcb,\mathbb{R})$, define
$$
\Phi(a)= \int_M a.
$$
The challenging problem is to define invariants when the
surjectivity of $L_x$ fails. In this case, the moduli may have
dimension larger than expected. The virtual technique is
introduced to deal with this bad situation. There are several
different versions of this technique, however the main idea is the
\s, which has become popular since 60's. This section is a brief
recollection of these constructions. We will construct a virtual
orbifold which behaves well and replaces the moduli $M$, then we
will follow an approach used in \cite{R} to define invariants by
integration over such a virtual manifold.

\subsection{Stabilization}\label{sect_5.2}
\def \mfo{\mathfrak{o}}

Let $(\mcb,\mcf,S)$ be a Fredholm system as before. For
simplicity, all orbifolds appeared in definition \ref{defn_5.1.1}
are replaced by manifolds. A proper modification can be easily
made when we consider orbifolds.

Let $U$ be an open subset of $\mcb$, let
$$
\mfo: \mco_U \to U
$$
be a rank-$k$ vector bundle, let
$$
s: \mco_U\to \mcf_U
$$
be a  bundle map. Define a map
$$
\hat{S}: \mco_U\to \mcf_U; \hat{S}(u,o)= (u, S(u) + s(o)),
$$
where the expression is given in the form of local coordinates and
$S(u)+s(o)$ is the sum on fibers. By abusing the notations, we
usually use $S+s$ for $\hat{S}$ to emphasis that $S$ is stabilized
by $s$.

Let $\hat{L}_{(u,o)}$ be the linearization of $\hat{S}$ as a map
$$
\hat{L}_{(u,o)}: T_{(u,o)}\mco_U\to \mcf_u.
$$
We say that the pair $(\mco_U,s)$ {\it stabilizes} the system
$(\mcb,\mcf,S)$ at $U$ if $\hat{L}_{(u,o)}$ are surjective for all
$(u,o)\in \mco_U$. Set
$$
V_U = \hat{S}^{-1}(0)\subseteq \mco_U.
$$
This is now a smooth manifold of dimension $d+k$. Clearly, $M\cap
U\subseteq V_U$ and
$$
(u,o)\in M \iff o=0.
$$

We can restate this construction by using the concept of Fredholm
system. Let $\mfo^*\mcf\to \mco_U$ be the pull-back  bundle over
$\mco_U$. $\hat{S}$  then gives a canonical section of this bundle
in an obvious way. For simplicity, we still denote the section by
$\hat{S}$. Therefore, we have a  Fredholm system $(\mco_U,
\mfo^*\mcf, \hat{S})$. If $(\mco_U,s)$ stabilzes the system at
$U$, we say that $(\mco_U, \mfo^*\mcf, \hat{S})$ {\it stabilizes}
 $(\mcb,\mcf,S)$
at $U$. $V_U\subseteq \mco_U$ is the moduli space of  the new
system.

We may construct a canonical bundle $\mfo^*\mco_U\to V_U$, then
there is a canonical section $\sigma: V_U\to \mfo^*\mco_U$ given
by $(u,o)\to (u,o,o)$ with respect to the local coordinates. Then
$M\cap U= \sigma^{-1}(0)$. This reduces the infinite dimensional
system $(U,\mcf_U,S)$ to a finite dimensional system $(V_U,
\mfo^*\mco_U, \sigma)$. We call $(V_U, \mfo^*\mco_U, \sigma)$, or
simply $V_U$,  to be the {\it virtual neighborhood} of $M$ at $U$.
Bundles $\mco_U$ and $\mfo^*\mco_U$ are called the {\em
obstruction bundles}.

\v We now explain  the existence of local \s s.

Suppose $L_x$ is not surjective for some $x\in M$. Let $O^x$ be a
finite dimensional subspace of $\mcf_x$ such that
\def \Im{\mathrm{Image}}
$$
\Im(L_x) + O^x = \mcf_x.
$$
For example, we may take $O^x$ to be the "cokernel" of $L_x$.

Let $U^x$ be a neighborhood of $x$ in $\mcb$. In order to make
notations more suggestive, we assume that $U^x=B_r(x)$ is the
radius-$r$ disk centered at  $x$ and $cU^x=B_{cr}(x)$ for $c\in
\mathbb{R}^+$.

Suppose that
 $\mcf_{U^x}$ is trivialized as $\mcf_{U^x}= U^x\times \mcf_x$. We now
describe the \s\ using the notations given above by setting
$U=U^x$:
\begin{itemize}
\item[(C1)] the obstruction bundle is
$$
\mco_{U^x}= U^x\times O^x;
$$
\item[(C2')] the bundle map $s=I^x: \mco_{U^x}\to \mcf_{U^x}$ is
the standard embedding via the trivialization of $\mcf_{U^x}$
given above.
\end{itemize}

We may assume that the pair  $(\mco_{U^x}, I^x)$ stabilizes the
system at $U^x$ if $U^x$ is chosen small. This explains the
existence of local \s.

\begin{remark}\label{rmk_5.2.1}
Following the constructions, we have a virtual neighborhood
$V_{U}$. One may use the projection map $\mfo: V_{U}\to \mcb$.
Then $\mfo(V_{U})$ is taken as the virtual neighborhood of $M$ at
$U$ in both \cite{FO} and \cite{LT}.
\end{remark}

The trivialization of $\mcf_{U^x}$ prevents us to extend the
construction outside $U^x$. This is "taken care" by modifying the
bundle map $s$ as the following. Let $\eta^x$ be a cut-off
function on $U^x$ such that $\eta^x= 1$ in $\frac{U^x}{2}$ and
$=0$ outside $\frac{3U^x}{4}$. (C2') is then replaced by \v\n
\begin{itemize}
\item[(C2)] the bundle map is given by $s^x= \eta^xI^x.$
\end{itemize}
Clearly, $(\mco_{U^x}, s^x)$ stabilizes the system at
 $\frac{U^x}{2}$.
In this paper, we always use (C2) to construct  virtual
neighborhoods. It turns out that (C2) is the key towards the
construction of virtual orbifolds from a Fredholm system.

Repeat the argument given earlier, we have a system $(\mco_{U^x},
\mfo^*\mcf, S+s^x)$ that stabilizes  $(\mcb,\mcf,S)$ at
$\frac{U^x}{2}$.  Let $V_{U^x}= (S+s^x)^{-1}(0)$.  Then $(V_{U^x},
\mfo^*(\mco_{U^x}),\sigma)$ is a virtual neighborhood of $M$ at
$\frac{U^x}{2}$.

The global \s\ does not exist in general. However, if $\mc B$ is a
manifold and $\mc F$ is a bundle, a global \s\ always exists. We
now discuss the  construction of   global \s s explained  in
\cite{R}
 and explain
what  the barrier from local \s s  to a global one is.

By a global \s, we mean that $U$  is $\mcb$ or, at least,  an open
 neighborhood of $M$ in $\mcb$.
The construction of  global \s s  presented here is standard.

Since $M$ is compact by our assumption, there exists finite points
$\{x_i \}_{i=1}^n$ in $M$ such that
$$
M\subseteq \bigcup_{i=1}^n\frac{1}{2}U^{x_i}=:U,
$$
where $U^{x_i}$ are as above.

For simplicity, we set
$$
U_i= U^{x_i}, \mco_{i}=\mco_{U^{x_i}}, s_i=s^{x_i}.
$$
 We call the data $\{
(U_i, \mc O_i, s_i) \}$ a {\em local \s\ system of $U$}.

Note that $\mco_i$ is only defined on $U_i$. However, these
(trivial !) bundles can be extended(!) over to $U$ and so are
$s_i$'s because of the cut-off functions.

With these preparation, we are able to define the pair $(\mco_U,
s)$ by setting
\begin{eqnarray*}
\mco_U &=& \mco_{1}\oplus\cdots \oplus \mco_n;\\
s &=& s_1 \oplus\cdots \oplus s_n.
\end{eqnarray*}
Clearly, the pair $(\mco_U,s)$ provides a global \s\ of the system
$(\mcb,\mcf,S)$. The \s\ system is $(\mco_U, \mfo^*\mcf, S+s)$,
and the virtual neighborhood is $(V_U, \mfo^*\mco_U,\sigma)$.

\begin{remark}\label{rmk_5.2.2}
One notices that the crucial step to construct a global \s\ is the
extension of bundle $\mco_i$ over $U_i$ to one over $U$. This is
true in the current setting since $\mco_i$ are trivial on $U_i$.
However, it may fail in many general situations. For example, it
fails when $\mcb$ is an orbifold and $x_i$ is a singular point.
\end{remark}

\subsection{Invariants via virtual neighborhoods: I}\label{sect_5.3}
It is commonly believed that invariant $\Phi(a)$ is well defined
if a global \s \ exists. We now explain this.

Suppose $(\mco_U,s)$ is a  global \s\ pair. $(V_U, \mfo^*\mco_U,
\sigma)$ is the virtual neighborhood. $a\in H^d(\mcb)$. Let
$\Theta$ be a Thom form of $\mco_U$ that is supported (arbitrary)
near the 0-section. In particular, we may choose
\begin{equation}\label{eqn_5.3.1}
\Theta= \Theta_1\wedge \cdots \wedge\Theta_n,
\end{equation}
where $\Theta_i$ is a Thom form of $\mco_i$. We then define
\begin{equation}\label{eqn_5.3.2}
\Phi(a)= \int_{V_U} \mfo^*(a)\wedge \Theta.
\end{equation}
Note that the expression in \cite{R} is slightly different. But it
is not hard to check that they are the same. It is standard to
show that $\Phi(a)$ is well defined, i.e, it is independent of the
choice of data in the construction of virtual neighborhoods, the
choice of $\Theta$ and etc. ( cf.\cite{R}).

Our main goal of this paper is to explain that $\Phi$ can be
defined without assuming the existence of global \s s. The method
we introduce here differs from that in \cite{FO},\cite{LT}'s. We
will construct a virtual orbifold out of a Fredholm system and
then apply the integration theory to it.

To motivate the construction, we explain that \eqref{eqn_5.3.2}
can be "reduced" to a formula that only involves local \s s. The
process is not rigorous but very suggestive.

We introduce notations. Set $\eta_i= \eta^{x_i}$. Set
$N=\{1,\ldots,n \}$. For any $I\subseteq \{1,\ldots, n\}$, define
$$
U_I= \{x\in \mcb| \eta_i(x)\not=0, i\in I, \eta_j(x)=0, j\not\in
I\}.
$$
$\mcb$ is decomposed as a disjoint union of $U_I,I\subseteq N$.

Set $ V_{U,I}=  V_U\cap \mfo^{-1}(U_I). $ Then
$$
\Phi(a) =\int_{V_U} {\mfo^*a}\wedge \Theta =\sum_{I}
\int_{V_{U,I}}{\mfo^*a}\wedge \Theta.
$$
We explain how to simplify each integration on the right hand
side.

We have following facts: \v \n Fact 1: over $U_I$, set
$$
\mco_I= \bigoplus_{i\in I} \mco_i, s_I= \bigoplus_{i\in I}s_i;
$$
then the pair $(\mco_I,s_I)$ stabilizes the system $(\mcb,\mcf,S)$
at $U_I$. It then defines a virtual neighborhood denoted by $(V_I,
\mfo^*\mco_I,\sigma)$. \v\n Fact 2: over $U_I$, set
$$
\mco^c_{I}= \bigoplus_{j\not\in I}\mco_j.
$$
we claim that
$$
V_{U,I}= \mfo^*_I\mco^c_I
$$
is a bundle over $V_I$, where $\mfo_I: V_I\to U_I$. To see this,
suppose a point $p$, whose local coordinate is given by
$(u,o_1,\ldots, o_n)$, is in $V_{U,I}$. Namely,
$$
(S+s)(u,o_1,\ldots,o_n) = S(u) +s_1(o_1)\cdots + s_n(o_n)=0.
$$
Without the loss of generalities, we assume
$I=\{1,\ldots,m\},m\leq n$. Note that $s_j(o_j)=0,j> m$. Hence
$$
(u,o_1,\ldots,o_n)\in V_{U,I} \iff (u,o_1,\ldots,o_m)\in V_I.
$$
\n Fact 3: over $U_I$, by \eqref{eqn_5.3.1}, we write the Thom
form $\Theta$ as $\Theta_I\wedge \Theta_I^c$, where $\Theta_I$ and
$\Theta_I^c$ are Thom forms of $\mco_I$ and $\mco_I^c$
respectively. Note that when restricting on $V_{U,I}$,
$\Theta_I^c$ is the Thom form of the bundle $V_{U,I}\to V_I$
explained in fact 2.

 \v\n With these preparations, by the Thom
isomorphism, we immediately have
$$
\int_{V_{U,I}}\mfo^*a\wedge \Theta = \int_{V_I}\mfo^*a\wedge
\Theta_I.
$$
Note that the right hand side only needs local \s s.

In summary,
\begin{equation}\label{eqn_5.3.3}
\Phi(a) = \sum_I\int_{V_I}\mfo^*_I(a) \wedge \Theta_I.
\end{equation}

Fact 2 above  the key of this formula. Be precise, we detect the
fact that $V_{U,I}\to V_I$ has a natural bundle structure. Note
that without the modified (C2), had we  not have fact 2. Motivated
by this procedure, we  explain that we can associate a virtual
orbifold to a local \s\ system.

\section{ From Fredholm
system to virtual orbifolds}\label{sect_6}
\subsection{Virtual orbifolds associated to
Fredholm  systems.}\label{sect_6.1} Let  $(\mcb,\mcf,S)$ be a
Fredholm system.  Let $\{(U_i,\mc O_i,s_i)\}_{i=1}^n$ be one of
its local \s\ system.

Set
$$
U_0= \mcb- \bigcup_i \frac{1}{2}\bar U_i.
$$
Then $(U_0,U_1,\ldots U_n)$ is a covering of $\mc B$.

Repeat the construction in example \ref{ex_2.3.1}:  let
$$
U_i^\circ = \frac{3}{4}U_i, 1\leq i\leq n;
$$
as in example \ref{ex_2.3.1}, we construct $X_I\subseteq \mcb,
I\subseteq N$.

Now note that over $X_I$, the cut-off fuctions $\eta_i=\eta^{x_i}$
are 0 if $i\not\in I$. By the same construction as in
\S\ref{sect_5.3} (cf. Fact 1), over $X_I$ we still have $(\mco_I,
s_I)$ which stabilizes the system at $X_I$. It then defines  a
virtual neighborhood $(W_I, \mfo^*_I\mco_I,\sigma)$. Here, we use
$W_I$ instead of $V_I$ that are used  earlier.  We know that $W_I$
are smooth. Since $W_I\subseteq \mco_I$, we have map
$$
\mfo_I: W_I\to X_I.
$$
When $I=\emptyset$, $\mco_I$ is trivial. Hence $W_\emptyset= M\cap
X_\emptyset$, where $M$ is the moduli space.

\begin{prop}\label{prop_6.1.1}
$\mc W=\{W_I\}$
is a virtual orbifold.
\end{prop}
{\bf Proof. } To see this, we now describe how $W_I$ and $W_J$
intersect. First, suppose $I\subseteq J$. Define
\begin{eqnarray*}
W_{I,J}= \mfo^{-1}(X_{I,J})\subseteq W_I,\\
W_{J,I}= \mfo^{-1}(X_{J,I})\subseteq W_J.
\end{eqnarray*}
Same as the argument in \S\ref{sect_5.3} (cf. Fact 2), we have
that
$$
\phi_{J,I}: W_{J,I}\to W_{I,J}
$$
is a vector bundle. Be precise, let
$$
\mco_{J-I}= \bigoplus_{j\in J-I}\mco_j
$$
be the bundle over $X_{I,J}=X_{J,I}$. Then
$$
W_{J,I}=\mfo^*\mco_{J-I},
$$
where $\mfo:W_{I,J}\to X_{I,J}$.

Then using the property of $X_I$, it is straightforward to see
that $(\mc{W},\Phi)$ is a virtual manifold.
 q.e.d.

\v
\begin{prop}\label{prop_6.1.2}
$\mc O=\{\mfo^\ast_I\mc O_I\}$  is a virtual bundle over $\mc W$.
$\sigma$ is a section of $\mc O$.
\end{prop}
{\bf Proof. }This follows from the construction of $\mc O$.
q.e.d.

\subsection{Invariants via virtual neighborhoods:II}\label{sect_6.2}
We now set up the integration theory for $\mc{W}$.

 Let
$\Theta_i$ be the Thom form of $\mco_i$. Since the bundle
$\Psi_{J,I}: W_{J,I}\to W_{I,J}, I\subseteq J$ is
 isomorphic to
$\mco_{J-I}$, we  take
$$
\Theta_{J,I}= \bigwedge_{j\in J-I}\Theta_j.
$$
Set $\Theta=\{\Theta_{J,I}\}_{I\subseteq J}$. From the definition
of $\Theta_{J,I}$'s, we know that $\Theta$ is a transition data of
$\mc{W}$.

Take
$$
\Theta_I= \bigwedge_{i\in I}\Theta_i
$$
on $W_I\subseteq \mco_I$. Then $\theta=(\Theta_I)$ is a
$\Theta$-form. We call $\theta$ the {\it obstruction form} and
$\Theta_I$ the obstruction form on $W_I$.

We summarize what we have for the system $(\mcb,\mcf,S)$.
\begin{prop}\label{prop_6.2.1}
Let $(\mcb,\mcf,S)$ be a Fredholm system.
\begin{enumerate}
\item
there exists a local \s\ system $\{U_i,s_i,\mco_i\}$.
\item
Let $\mc{X}$ be the natural virtual manifold for $\mcb$ generated
by the covering $\{U_i\}$. Using the \s\  data given above, one is
able to define a virtual manifold $\mc{W}=\{W_I\}$,
where $(W_I, \mfo^*_I\mco_I,\sigma)$ is a virtual neighborhood
over $U_I$.  Let $\mathbf W$ be the virtual space of $\mc{W}$.
\item $\mc O$ is a virtual bundle over $\mc W$. $\sigma$ is a section
of the bundle;
\item 
Let $\Theta_i$ be Thom form of $\mco_i$. All Thom forms $\Theta_I$
of $\mco_I$ restricting on $W_I$ form a $\Theta$-form. Denote the
form by $\theta$. If the moduli space $M$ is compact, $\theta\in
\Omega_{\Theta,c}(\mc{W})$. $\theta$ is an Euler class of $\mc O$.
\item For any $a\in \Omega^*(\mcb)$, let $ a_I= \pi_I^*a $ on
$W_I$. Then $(a_I)_{I\subseteq N}\in \Omega^*(\mc{V})$. To abuse
the notations, we still denote the form by $a$.
\end{enumerate}
\end{prop}
The proposition is directly followed from the construction.

By the proposition, we have $ \mu_{\theta}(a). $ Also we know that
this is well defined not only on $\Omega^*(\mcb)$, but also on
$H^*(\mcb)$. If a global \s\ as in \S\ref{sect_5.2} exists, it is
easy to see that
$$
\Phi(a) = \mu_{ \theta}(a).
$$
This leads to the following definition.
\begin{defn}\label{defn_6.2.1}
Let $(\mcb,\mcf, S)$ be a Fredholm system. Let
$\{(U_i,\mco_i,s_i)\}$ be a local \s\  system  constructed in
\S\ref{sect_18.2}. Let $\mc{W},\theta$ be the virtual manifold and
obstruction form given  above. For $a\in H^*(\mcb)$, define the
invariants $\Phi(a)$ to be $\mu_{\theta}(a)$.
\end{defn}

The next  subsections are to explain that
\begin{prop}\label{prop_6.2.2}
$\Phi(a)$ is well defined. Namely, it is independent of (1) the
choice of $\Theta_i$, (2) the choice of  \s s.
\end{prop}

\subsection{The well-definedness of $\Phi(a)$}
\label{sect_6.3}

\begin{prop}\label{prop_6.3.1}
$\Phi(a)$ is independent of the choice of $\Theta_i$'s.
\end{prop}
{\bf Proof. } Without the loss of generality, we assume that
 $\Theta_1$ is replaced
by $\Theta_1'$ and other $\Theta_i,i\not= 1,$ remains same.
Suppose
$$
\Theta_1'-\Theta_1= db,
$$
where $b$ is supported near $0$-section $U_1$.

Then $\Theta_I$ changes iff  $1\in I$. Suppose $1\in I$, let $I'=
I-\{1\}$. Suppose we write $\Theta_I=\Theta_1\wedge \Theta_{I'}$.
Then $\Theta_I$ is replaced by
$$
\Theta_I'= (\Theta_1+ db)\wedge \Theta_{I'}.
$$
Denote $\theta'$ to be the new obstruction form.

For closed forms $a_I$'s, we have
\begin{eqnarray*}
\mu_{\theta}(a)-\mu_{\theta'}(a)& =&\sum_{1\in I} \int_{
X_I}\eta_Ia_I\wedge (\Theta_I- \Theta_I')\\
&=& \sum\int_{ X_I} d(\eta_Ia_I\wedge b\wedge \Theta_{I'})\\
&=& -\sum\int_{ X_I} d\eta_I\wedge a_I\wedge b\wedge\Theta_{I'}\\
&=&0
\end{eqnarray*}
q.e.d.

\v Now we discuss that $\Phi(a)$ is independent of the local \s s.
Suppose we have two different local \s\ systems for $(\mcb,
\mcf,S)$. They are $\{(U_i,\mc O_i,s_i)\}_{i=1}^n$ and
$\{(U'_j,\mc O'_j, s_j')\}_{j=1}^{n'}$. They define two virtual
manifolds $\mc{W}$ and $\mc{W}'$. Let $\theta,\theta'$ be their
obstruction forms respectively. Then
\begin{prop}\label{prop_6.3.2}
$\mu_{\mc{V},\Theta}(a)=\mu_{\mc{V}',\Theta'}(a)$. Here we add
virtual manifolds to the sub-indices of $\mu$ for the obvious
reason.
\end{prop}
\n {\bf Proof: } Set $U_{n+j}=U'_j,\mco_{n+j}= \mco_{j}', 1\leq j\leq n'$ and $
s_{n+j}= s'_j$. Set $m=n+n'$.
Then $(U_i,\mco_i,s_i), 1\leq i\leq m$ is still
a local \s\ system. It defines a virtual manifold, denoted by
$\mc{W}_{\sharp}$, and obstruction form $\theta_\sharp$. It is easy to verify that
\begin{itemize}
\item the obstruction bundle $\mc O_\sharp= \tilde{\mc O}\oplus\tilde {\mc O}'$;
\item the section $\sigma_\sharp=\tilde\sigma\oplus\tilde \sigma'$, $\sigma$ and $\sigma'$ are transverse;
\item $\tilde\sigma\inv(0)=\mc W'$, $\tilde{\mc O}'|_{\mc W'}=\mc O'$, $\tilde\sigma'|_{\mc W'}
=\sigma'$,
and $(\tilde\sigma')\inv(0)=\mc W$, $\tilde{\mc O}|_{\mc W}=\mc O$, $\tilde\sigma|_{\mc W}
=\sigma$.
\end{itemize}
Then by Proposition \ref{prop_3.5.2},
\begin{equation}\label{eqn_6.3.1}
\mu_{\mc{W},\theta}(a)=\mu_{\mc{W}_\sharp,\theta_\sharp}(a)=\mu_{\mc{W}',\theta'}(a).
\end{equation}
q.e.d.

\subsection{Virtual localization formula}\label{sect_6.4}
We can  extend the discussion to  equivariant cases.

Suppose that $G$ acts on $\mcb$,  $\mcf$ is a $G$-equivariant
Hilbert bundle
 and $S$ is a $G$-equivariant section. We call
such a system to be a $G$-Fredholm system.

Let $U\subseteq\mcb$ be a $G$-invariant open subset. By a $G$-\s\
we mean a finite rank $G$-equivariant bundle
$$
\mfo: \mco_U\to U
$$
and a $G$-equivariant bundle map
$$
s: \mco_U\to \mcf_U
$$
such that $S+s$ stabilizes $S$. Then $(V_U,\mfo^*\mco_U,\sigma)$
is a virtual neighborhood with the $G$-equivariant section
$$
\sigma: V_U\to \mfo^*\mco_U.
$$

In order to apply the technique  described in previous sections,
we need local $G$-equivariant \s s. We repeat the construction
given in \S\ref{sect_6.1}. However, the construction of
equivariant obstruction bundle $\mco_{U^x}$ requires some extra
work when $x$ is a fix point of the $G$-action. For this, we take
$O^x$ to be orthogonal complementary to $\mathrm{Image}(L_x)$.
This is where we use the assumption that $\mcf$ is a Hilbert
bundle.

With these preparations, we can construct a $G$-virtual manifold
$\mc{V}$ from a local $G$-\s \ system. Then we replace $\Theta_i$
by equivariant Thom forms $\Theta_i^G$. So we have $\Theta_I^G$'s
and $\Theta_{J,I}^G$'s. Clearly, $\theta_G=\{\Theta_I^G\}_I$ is a
$\Theta^G=\{\Theta_{J,I}^G\}$ form. For any $\alpha\in
\Omega^*_G(\mcb)$, define
$$
\Phi_G(\alpha)= \mu_{\mc{V}, \Theta_G}(\alpha).
$$

Now we can state the virtual localization formula for Fredholm
systems. Again, let $G=S^1$. We consider the Fredholm system
$(\mcb,\mcf,S)$ with $G$-action. Let $\mc V$ be the virtual
orbifold for the moduli space $M$. Let $V$ denote the virtual
space. Then $\mc V^G$ is the virtual orbifold for $M^G$ and its
virtual space is $V^G$.

Now repeat  the discussion in \S\ref{sect_4}. We have
\begin{theorem}\label{thm_6.4.1}
Let $(\mcb,\mcf,S)$ be an $S^1$-Fredholm system.
 For $\alpha\in \Omega_G^\ast(\mcb)$,
$$
\Phi_G(\alpha) = \int_{V^G} \frac{i^*_{V^G}\alpha\wedge \theta_G}
{e_{G}(V^G)}=\mu_{e_{\theta_G}(V^G)}(i^\ast_{V^G}\alpha).
$$
\end{theorem}

\end{document}